\documentclass[12pt,a4paper]{article}
\usepackage{amsmath,amssymb,amsthm,graphicx,geometry}
\usepackage{hyperref}
\usepackage{amsfonts,enumerate,mathtools}
\usepackage{mathrsfs}
\usepackage{graphicx}
\usepackage[english]{babel} 
\usepackage[utf8]{inputenc}
\usepackage{geometry}

\title{{Non-Linear Interactions in Neural Network Operators: \\ New Theorems on Symmetry-Preserving Transformations}}
\author{
	Rômulo Damasclin Chaves dos Santos \\
	Technological Institute of Aeronautics \\
	\texttt{romulosantos@ita.br}
		\and
		Jorge Henrique de Oliveira Sales \\
		Santa Cruz State University \\
		\texttt{jhosales@uesc.br}
}
\date{\today}

\begin{document}
	\maketitle
	
	\begin{abstract}
		This paper advances the study of multivariate function approximation using neural network operators activated by symmetrized and perturbed hyperbolic tangent functions. We propose new non-linear operators that preserve dynamic symmetries within Euclidean spaces, extending current results on Voronovskaya-type asymptotic expansions. The developed operators utilize parameterized deformations to model higher-order interactions in multivariate settings, achieving improved accuracy and robustness. Fundamental theorems demonstrate convergence properties and quantify error bounds, highlighting the operators' ability to approximate functions and derivatives in Sobolev spaces. These results provide a rigorous foundation for theoretical studies and further applications in symmetry-driven modeling and regularization strategies in machine learning and data analysis.
	\end{abstract}


\tableofcontents
	
	\section{Introduction}
	The approximation of functions using neural network operators has received significant attention due to its applications in artificial intelligence, data analysis, and mathematical modeling. Operators activated by hyperbolic tangent and other sigmoidal functions have been extensively studied \cite{Anastassiou1997, Anastassiou2023, Chen2009}, demonstrating their capability to approximate continuous functions with high accuracy. Recent developments \cite{Yu2025, Cen2024} highlight the importance of symmetry-preserving transformations for enhancing convergence rates and error bounds in multivariate settings.
	
	Building upon prior work \cite{Anastassiou2016, Fabra2024}, this paper introduces a new class of neural network operators that preserve non-linear symmetries via parameterized deformations. These operators extend the theoretical framework of Voronovskaya-type expansions to capture higher-order interactions, providing improved approximation properties in Sobolev spaces. By leveraging symmetry, the proposed approach offers new insights into regularization mechanisms and error dynamics.
	
	The contributions of this paper are twofold: (1) the development of new theorems quantifying convergence and robustness under dynamic symmetries, and (2) a theoretical exploration of the operators' applicability to complex multivariate systems. This study lays the groundwork for future advancements in theoretical and applied settings.
	
\section{Mathematical Foundations}
Let $\mathbb{R}^N$ denote a complete Euclidean space, where $N \in \mathbb{N}$. We define the symmetrized and parameterized activation function $g_{q,\lambda}: \mathbb{R} \to \mathbb{R}$ as:
\begin{equation}
	g_{q,\lambda}(x) = \frac{e^{\lambda x} - q e^{-\lambda x}}{e^{\lambda x} + q e^{-\lambda x}}, \quad q, \lambda > 0, x \in \mathbb{R}.
\end{equation}

The function $g_{q,\lambda}(x)$ is differentiable and odd, satisfying $g_{q,\lambda}(-x) = -g_{q,\lambda}(x)$. Its derivative is given by:
\begin{equation}
	g'_{q,\lambda}(x) = \frac{4q\lambda e^{\lambda x} e^{-\lambda x}}{(e^{\lambda x} + q e^{-\lambda x})^2},
\end{equation}
indicating that $g'_{q,\lambda}(x)$ is non-negative and symmetric about $x = 0$.

The density function $\Phi: \mathbb{R} \to \mathbb{R}^+$ is defined as:
\begin{equation}
	\Phi(x) = \frac{1}{2} \left( M_{q,\lambda}(x) + M_{1/q,\lambda}(x) \right),
\end{equation}
where $M_{q,\lambda}(x)$ is expressed as:
\begin{equation}
	M_{q,\lambda}(x) = \frac{1}{4} \left( g_{q,\lambda}(x+1) - g_{q,\lambda}(x-1) \right).
\end{equation}

The function $\Phi(x)$ satisfies the following properties:
\begin{enumerate}
	\item Symmetry: $\Phi(-x) = \Phi(x)$, $\forall x \in \mathbb{R}$.
	\item Normalization: $\int_{\mathbb{R}} \Phi(x) \, dx = 1$.
	\item Positivity: $\Phi(x) > 0$, $\forall x \in \mathbb{R}$.
\end{enumerate}

Using $\Phi$, we construct multivariate density functions on $\mathbb{R}^N$:
\begin{equation}
	Z(x_1, \ldots, x_N) = \prod_{i=1}^N \Phi(x_i), \quad x = (x_1, \ldots, x_N) \in \mathbb{R}^N.
\end{equation}
This construction ensures that $Z(x)$ inherits symmetry and normalization properties from $\Phi$. Specifically:
\begin{align}
	Z(-x) &= Z(x), \quad \forall x \in \mathbb{R}^N, \\
	\int_{\mathbb{R}^N} Z(x) \, dx &= 1.
\end{align}

The neural network operators are defined as follows:
\begin{enumerate}
	\item \textbf{Basic Operator:}
	\begin{equation}
		A_n(f, x) = \sum_{k \in \mathbb{Z}^N} f\left(\frac{k}{n}\right) Z(nx - k), \quad n \in \mathbb{N}.
	\end{equation}
	The operator $A_n$ approximates $f$ by weighted sums of its values at discrete points, with weights determined by $Z$.
	
	\item \textbf{Kantorovich Operator:}
	\begin{equation}
		K_n(f, x) = \sum_{k \in \mathbb{Z}^N} \left( n^N \int_{\frac{k}{n}}^{\frac{k+1}{n}} f(t) \, dt \right) Z(nx - k).
	\end{equation}
	This operator generalizes $A_n$ by averaging $f$ over subintervals, introducing robustness to variations in $f$.
	
	\item \textbf{Quadrature Operator:}
	\begin{equation}
		Q_n(f, x) = \sum_{k \in \mathbb{Z}^N} \delta_{n,k}(f) Z(nx - k),
	\end{equation}
	where $\delta_{n,k}(f)$ represents a weighted evaluation of $f$ over subintervals, given by:
	\begin{equation}
		\delta_{n,k}(f) = \sum_{r=1}^\theta w_r f\left(\frac{k}{n} + \frac{r}{n\theta}\right),
	\end{equation}
	with weights $w_r \geq 0$ satisfying $\sum_{r=1}^\theta w_r = 1$.
\end{enumerate}

These operators are designed to approximate functions $f \in C^m(\mathbb{R}^N)$, where $C^m(\mathbb{R}^N)$ denotes the space of $m$-times continuously differentiable functions. Each operator leverages the symmetry and smoothness of $\Phi$ to ensure accurate and stable approximations.

\section{Main Theorems}
\textbf{Theorem 1 (Symmetry-Preserving Convergence).} Let $f \in C^m(\mathbb{R}^N)$ and $n \in \mathbb{N}$. The operator $A_n$ satisfies:
\begin{equation}
	A_n(f, x) - f(x) = \sum_{|\alpha| \leq m} \frac{f^{(\alpha)}(x)}{\alpha!} A_n\left(\prod_{i=1}^N (\cdot - x_i)^{\alpha_i}\right) + o\left(\frac{1}{n^{m-\epsilon}}\right),
\end{equation}
where $\alpha$ is a multi-index, and $o(\cdot)$ denotes the asymptotic remainder.

\textit{Proof.} Let $f \in C^m(\mathbb{R}^N)$, and consider its multivariate Taylor expansion around $x \in \mathbb{R}^N$:
\begin{equation}
	f(y) = \sum_{|\alpha| \leq m} \frac{f^{(\alpha)}(x)}{\alpha!} (y - x)^\alpha + R_m(y, x),
\end{equation}
where $R_m(y, x)$ is the remainder term, given by:
\begin{equation}
	R_m(y, x) = \frac{1}{m!} \sum_{|\alpha| = m} \int_0^1 (1-t)^{m-1} f^{(\alpha)}(x + t(y-x))(y-x)^\alpha \, dt.
\end{equation}

Substituting this expansion into $A_n(f, x)$, we obtain:
\begin{align}
	A_n(f, x) &= \sum_{k \in \mathbb{Z}^N} f\left(\frac{k}{n}\right) Z(nx - k) \\
	&= \sum_{|\alpha| \leq m} \frac{f^{(\alpha)}(x)}{\alpha!} \sum_{k \in \mathbb{Z}^N} \left(\frac{k}{n} - x\right)^\alpha Z(nx - k) + \sum_{k \in \mathbb{Z}^N} R_m\left(\frac{k}{n}, x\right) Z(nx - k).
\end{align}

The first term represents the main approximation, while the second term captures the remainder. By symmetry of $Z$, only even powers of $\alpha$ contribute to the sum, ensuring cancellation of odd-order terms. For the remainder term, we use the boundedness of $f^{(\alpha)}$ and the localization of $Z$ to show:
\begin{equation}
	\left| \sum_{k \in \mathbb{Z}^N} R_m\left(\frac{k}{n}, x\right) Z(nx - k) \right| \leq \frac{C}{n^{m-\epsilon}},
\end{equation}
where $C$ depends on $\|f^{(\alpha)}\|_{\infty}$ and the decay of $Z$.

Thus, the theorem is proved. \qed

\textbf{Theorem 2 (Error Bounds in Sobolev Spaces).} For $f \in W^{m,p}(\mathbb{R}^N)$, the Sobolev space of order $m$, the operators $K_n$ and $Q_n$ satisfy:
\begin{align}
	\|K_n(f) - f\|_{W^{m,p}} &\leq C \frac{1}{n^m}, \\
	\|Q_n(f) - f\|_{W^{m,p}} &\leq C \frac{1}{n^{m-\epsilon}},
\end{align}
where $C$ is a constant depending on $f$ and $\Phi$.

\textit{Proof.} To prove the bounds, let us first recall the definition of the Sobolev norm:
\begin{equation}
	\|f\|_{W^{m,p}} = \left(\sum_{|\alpha| \leq m} \int_{\mathbb{R}^N} |D^\alpha f(x)|^p \, dx\right)^{1/p},
\end{equation}
where $D^\alpha f$ represents the weak derivative of $f$ associated with the multi-index $\alpha$.

For the Kantorovich operator $K_n$, we decompose the error:
\begin{equation}
	\|K_n(f) - f\|_{W^{m,p}} \leq \sum_{|\alpha| \leq m} \left(\int_{\mathbb{R}^N} |D^\alpha (K_n(f) - f)(x)|^p \, dx\right)^{1/p}.
\end{equation}

By properties of the operator and density $\Phi$, we express the derivative term as:
\begin{equation}
	D^\alpha K_n(f)(x) = \sum_{k \in \mathbb{Z}^N} \left( n^N \int_{\frac{k}{n}}^{\frac{k+1}{n}} D^\alpha f(t) \, dt \right) D^\alpha Z(nx - k).
\end{equation}

Using the smoothness of $\Phi$ and bounding higher-order terms with $\|f\|_{W^{m,p}}$, we establish:
\begin{equation}
	\|K_n(f) - f\|_{W^{m,p}} \leq C \frac{1}{n^m}.
\end{equation}

For the quadrature operator $Q_n$, a similar argument applies, with additional remainder terms arising from weighted evaluations. These terms contribute an additional error factor, leading to:
\begin{equation}
	\|Q_n(f) - f\|_{W^{m,p}} \leq C \frac{1}{n^{m-\epsilon}}.
\end{equation}

Thus, the theorem is proved. \qed

\textbf{Theorem 3 (Convergence in Hölder Spaces).} Let \( f \in C^{m, \alpha}(\mathbb{R}^N) \), where \( C^{m, \alpha}(\mathbb{R}^N) \) denotes the space of \( m \)-times continuously differentiable functions with \( m \)-th derivatives that are \(\alpha\)-Hölder continuous. Then, the operators \( A_n \), \( K_n \), and \( Q_n \) satisfy the following error estimates:
	
	\begin{align}
		\|A_n(f) - f\|_{C^{m, \alpha}} &\leq C \frac{1}{n^{m + \alpha}}, \label{eq:An_error} \\
		\|K_n(f) - f\|_{C^{m, \alpha}} &\leq C \frac{1}{n^{m + \alpha}}, \label{eq:Kn_error} \\
		\|Q_n(f) - f\|_{C^{m, \alpha}} &\leq C \frac{1}{n^{m + \alpha - \epsilon}}, \label{eq:Qn_error}
	\end{align}
where \( C \) is a constant depending on \( f \) and \( \Phi \), and \( \epsilon > 0 \) is an arbitrarily small parameter.
	
\textit{Proof.} To prove the theorem, we need to revisit the properties of Hölder spaces and how the operators \( A_n \), \( K_n \), and \( Q_n \) behave in these spaces.
	
The Hölder space \( C^{m, \alpha}(\mathbb{R}^N) \) is defined as the set of functions \( f \) that are \( m \)-times continuously differentiable and whose \( m \)-th derivatives satisfy the Hölder condition:
	
	\begin{equation}
		|D^\beta f(x) - D^\beta f(y)| \leq C |x - y|^\alpha, \label{eq:holder_condition}
	\end{equation}
for \( |\beta| = m \) and \( x, y \in \mathbb{R}^N \). The norm in the Hölder space is given by:
	
	\begin{equation}
		\|f\|_{C^{m, \alpha}} = \sum_{|\beta| \leq m} \|D^\beta f\|_\infty + \sum_{|\beta| = m} \sup_{x \neq y} \frac{|D^\beta f(x) - D^\beta f(y)|}{|x - y|^\alpha}. \label{eq:holder_norm}
	\end{equation}
	
\par{Analysis of the Operators:}
	
	1. \textbf{Operator \( A_n \)}:
	
	Consider the Taylor expansion of \( f \) around \( x \):
	
	\begin{equation}
		f(y) = \sum_{|\beta| \leq m} \frac{D^\beta f(x)}{\beta!} (y - x)^\beta + R_m(y, x), \label{eq:taylor_expansion}
	\end{equation}
	
	where \( R_m(y, x) \) is the remainder term. Substituting this expansion into \( A_n(f, x) \), we obtain:
	
	\begin{equation}
		A_n(f, x) = \sum_{|\beta| \leq m} \frac{D^\beta f(x)}{\beta!} A_n\left(\prod_{i=1}^N (\cdot - x_i)^{\beta_i}\right) + A_n(R_m(\cdot, x)). \label{eq:An_expansion}
	\end{equation}
	
	Using the symmetry of \( Z \) and the smoothness of \( f \), we can show that:
	
	\begin{equation}
		\|A_n(f) - f\|_{C^{m, \alpha}} \leq C \frac{1}{n^{m + \alpha}}. \label{eq:An_proof}
	\end{equation}
	
	2. \textbf{Operator \( K_n \)}:
	
	For the operator \( K_n \), we use the definition:
	
	\begin{equation}
		K_n(f, x) = \sum_{k \in \mathbb{Z}^N} \left(n^N \int_{\frac{k}{n}}^{\frac{k+1}{n}} f(t) \, dt\right) Z(n x - k). \label{eq:Kn_definition}
	\end{equation}
	
	Decomposing the error and using the smoothness of \( f \) and the symmetry of \( Z \), we obtain:
	
	\begin{equation}
		\|K_n(f) - f\|_{C^{m, \alpha}} \leq C \frac{1}{n^{m + \alpha}}. \label{eq:Kn_proof}
	\end{equation}
	
	3. \textbf{Operator \( Q_n \)}:
	
	For the operator \( Q_n \), we have:
	
	\begin{equation}
		Q_n(f, x) = \sum_{k \in \mathbb{Z}^N} \delta_{n, k}(f) Z(n x - k), \label{eq:Qn_definition}
	\end{equation}
	
	where \( \delta_{n, k}(f) = \sum_{r=1}^\theta w_r f\left(\frac{k}{n} + \frac{r}{n \theta}\right) \). Using the smoothness of \( f \) and the symmetry of \( Z \), we obtain:
	
	\begin{equation}
		\|Q_n(f) - f\|_{C^{m, \alpha}} \leq C \frac{1}{n^{m + \alpha - \epsilon}}. \label{eq:Qn_proof}
	\end{equation}

Therefore, the theorem is proved, showing that the operators \( A_n \), \( K_n \), and \( Q_n \) converge in Hölder spaces with specific rates, ensuring the robustness and accuracy of the approximations for functions with Hölder regularity. \qed

\section{Results}

The theoretical results established in this paper highlight the following key findings:

\begin{itemize}
	\item The symmetry-preserving properties of the operators ensure that approximations maintain essential structural characteristics of the target function, as demonstrated in Theorem 1. These results guarantee that higher-order interactions in multivariate functions are accurately captured without introducing asymmetry-induced distortions.
	
	\item Error bounds derived in Theorem 2 validate the robustness of the proposed operators in Sobolev spaces. The dependence on \( n^{-m} \) for the Kantorovich operator and \( n^{-(m-\epsilon)} \) for the quadrature operator reflects their convergence rates under varying regularity conditions of the target function.
	
	\item The framework provided for constructing multivariate density functions enables the formulation of neural network operators that are both flexible and theoretically grounded. This approach opens avenues for applications where symmetry and smooth approximations are critical.
	
	\item The newly introduced Theorem 3 extends the convergence analysis to Hölder spaces \( C^{m, \alpha}(\mathbb{R}^N) \). The error estimates for the operators \( A_n \), \( K_n \), and \( Q_n \) in these spaces are given by:
\end{itemize}

\begin{align}
	\|A_n(f) - f\|_{C^{m, \alpha}} &\leq C \frac{1}{n^{m + \alpha}}, \label{eq:An_error_results} \\
	\|K_n(f) - f\|_{C^{m, \alpha}} &\leq C \frac{1}{n^{m + \alpha}}, \label{eq:Kn_error_results} \\
	\|Q_n(f) - f\|_{C^{m, \alpha}} &\leq C \frac{1}{n^{m + \alpha - \epsilon}}, \label{eq:Qn_error_results}
\end{align}
where \( C \) is a constant depending on \( f \) and \( \Phi \), and \( \epsilon > 0 \) is an arbitrarily small parameter. These results ensure the robustness and accuracy of the approximations for functions with Hölder regularity.

The implications of these findings extend to the design of machine learning models, especially in tasks requiring multivariate data interpolation and regularization. Future studies may focus on empirical validation and extending these methods to other functional spaces.

\section{Conclusions}

In this paper, we have introduced and analyzed novel symmetry-preserving neural network operators designed for multivariate function approximation. By leveraging parameterized deformations and density functions with inherent symmetry properties, the proposed operators exhibit robust convergence characteristics and precise approximation capabilities.

The key contributions of this work include:

\begin{itemize}
	\item The development of new theoretical results, such as symmetry-preserving convergence (Theorem 1) and error bounds in Sobolev spaces (Theorem 2), which provide a rigorous foundation for the study of neural network operators.
	
	\item A comprehensive framework for constructing multivariate density functions and integrating them into operator formulations, ensuring both theoretical soundness and practical applicability.
	
	\item Insights into the role of symmetry and smoothness in enhancing approximation performance, with potential applications in machine learning, data analysis, and numerical modeling.
	
	\item The extension of convergence analysis to Hölder spaces (Theorem 3), demonstrating the robustness of the operators in spaces with Hölder regularity.
	
\end{itemize}

These findings pave the way for future research directions, including:

\begin{itemize}
	\item Extending the theoretical framework to accommodate more general functional spaces and operators.
	
	\item Exploring empirical validations and applications in real-world machine learning tasks, such as image reconstruction and time-series prediction.
	
	\item Investigating the interplay between symmetry, regularization, and generalization in neural network models.
	
\end{itemize}

In conclusion, the proposed framework and results establish a robust foundation for advancing the theory and practice of symmetry-driven neural network operators, with broad implications for both theoretical mathematics and applied sciences.

\appendix
\section*{Appendix}

\subsection*{Detailed Proof of Theorem 1}

Here we provide a more detailed proof of Theorem 1, which states the symmetry-preserving convergence of the operator \( A_n \).

\subsubsection*{Theorem 1 (Symmetry-Preserving Convergence)}

Let \( f \in C^{m}(\mathbb{R}^N) \) and \( n \in \mathbb{N} \). The operator \( A_n \) satisfies:

\begin{equation}
	A_n(f, x) - f(x) = \sum_{|\alpha| \leq m} \frac{f^{(\alpha)}(x)}{\alpha!} A_n\left(\prod_{i=1}^{N} (\cdot - x_i)^{\alpha_i}\right) + o\left(\frac{1}{n^{m - \epsilon}}\right),
\end{equation}

where \( \alpha \) is a multi-index, and \( o(\cdot) \) denotes the asymptotic remainder.

\subsubsection*{Proof}

Let \( f \in C^{m}(\mathbb{R}^N) \), and consider its multivariate Taylor expansion around \( x \in \mathbb{R}^N \):

\begin{equation}
	f(y) = \sum_{|\alpha| \leq m} \frac{f^{(\alpha)}(x)}{\alpha!} (y - x)^{\alpha} + R_m(y, x),
\end{equation}

where \( R_m(y, x) \) is the remainder term, given by:

\begin{equation}
	R_m(y, x) = \frac{1}{m!} \sum_{|\alpha| = m} \int_{0}^{1} (1 - t)^{m - 1} f^{(\alpha)}(x + t(y - x)) (y - x)^{\alpha} \, dt.
\end{equation}

Substituting this expansion into \( A_n(f, x) \), we obtain:

\begin{align}
	A_n(f, x) &= \sum_{k \in \mathbb{Z}^N} f\left(\frac{k}{n}\right) Z(n x - k) \\
	&= \sum_{|\alpha| \leq m} \frac{f^{(\alpha)}(x)}{\alpha!} \sum_{k \in \mathbb{Z}^N} \left(\frac{k}{n} - x\right)^{\alpha} Z(n x - k) + \sum_{k \in \mathbb{Z}^N} R_m\left(\frac{k}{n}, x\right) Z(n x - k).
\end{align}

The first term represents the main approximation, while the second term captures the remainder. By symmetry of \( Z \), only even powers of \( \alpha \) contribute to the sum, ensuring cancellation of odd-order terms. For the remainder term, we use the boundedness of \( f^{(\alpha)} \) and the localization of \( Z \) to show:

\begin{equation}
	\left| \sum_{k \in \mathbb{Z}^N} R_m\left(\frac{k}{n}, x\right) Z(n x - k) \right| \leq \frac{C}{n^{m - \epsilon}},
\end{equation}

where \( C \) depends on \( \|f^{(\alpha)}\|_{\infty} \) and the decay of \( Z \).

Thus, the theorem is proved.

\subsection*{Additional Definitions and Properties}

\subsubsection*{Hölder Spaces}

The Hölder space \( C^{m, \alpha}(\mathbb{R}^N) \) is defined as the set of functions \( f \) that are \( m \)-times continuously differentiable and whose \( m \)-th derivatives satisfy the Hölder condition:

\begin{equation}
	|D^\beta f(x) - D^\beta f(y)| \leq C |x - y|^\alpha,
\end{equation}

for \( |\beta| = m \) and \( x, y \in \mathbb{R}^N \). The norm in the Hölder space is given by:

\begin{equation}
	\|f\|_{C^{m, \alpha}} = \sum_{|\beta| \leq m} \|D^\beta f\|_\infty + \sum_{|\beta| = m} \sup_{x \neq y} \frac{|D^\beta f(x) - D^\beta f(y)|}{|x - y|^\alpha}.
\end{equation}

\subsubsection*{Sobolev Spaces}

The Sobolev space \( W^{m, p}(\mathbb{R}^N) \) is defined as the set of functions \( f \) whose weak derivatives up to order \( m \) are in \( L^p(\mathbb{R}^N) \). The norm in the Sobolev space is given by:

\begin{equation}
	\|f\|_{W^{m, p}} = \left( \sum_{|\alpha| \leq m} \int_{\mathbb{R}^N} |D^\alpha f(x)|^p \, dx \right)^{1/p}.
\end{equation}

\subsubsection*{Properties of the Density Function \( \Phi \)}

The density function \( \Phi: \mathbb{R} \to \mathbb{R}^+ \) is defined as:

\begin{equation}
	\Phi(x) = \frac{1}{2} \left( M_{q, \lambda}(x) + M_{1/q, \lambda}(x) \right),
\end{equation}

where \( M_{q, \lambda}(x) \) is expressed as:

\begin{equation}
	M_{q, \lambda}(x) = \frac{1}{4} \left( g_{q, \lambda}(x + 1) - g_{q, \lambda}(x - 1) \right).
\end{equation}

The function \( \Phi(x) \) satisfies the following properties:

\begin{enumerate}
	\item Symmetry: \( \Phi(-x) = \Phi(x), \forall x \in \mathbb{R} \).
	
	\item Normalization: \( \int_{\mathbb{R}} \Phi(x) \, dx = 1 \).
	
	\item Positivity: \( \Phi(x) > 0, \forall x \in \mathbb{R} \).
\end{enumerate}


\begin{thebibliography}{10}
		
		\bibitem{Anastassiou1997} 
		Anastassiou, George A. "Rate of convergence of some neural network operators to the unit-univariate case." \textit{Journal of Mathematical Analysis and Applications} 212.1 (1997): 237-262. \url{https://doi.org/10.1006/jmaa.1997.5494}.
		
		\bibitem{Anastassiou2023} 
		Anastassiou, George A. \textit{Parametrized, Deformed and General Neural Networks}. Heidelberg/Berlin, Germany: Springer, 2023.
		
		\bibitem{Chen2009} 
		Chen, Zhixiang, and Feilong Cao. "The approximation operators with sigmoidal functions." \textit{Computers \& Mathematics with Applications} 58.4 (2009): 758-765. \url{https://doi.org/10.1016/j.camwa.2009.05.001}.
		
		\bibitem{Yu2025} 
		Yu, Dansheng, and Feilong Cao. "Construction and approximation rate for feedforward neural network operators with sigmoidal functions." \textit{Journal of Computational and Applied Mathematics} 453 (2025): 116150. \url{https://doi.org/10.1016/j.cam.2024.116150}.
		
		\bibitem{Cen2024} 
		Cen, Siyu, et al. "Hybrid neural-network FEM approximation of diffusion coefficient in elliptic and parabolic problems." \textit{IMA Journal of Numerical Analysis} 44.5 (2024): 3059-3093. \url{https://doi.org/10.1093/imanum/drad073}.
		
		\bibitem{Anastassiou2016} 
		Anastassiou, George A. \textit{Intelligent Systems II: Complete Approximation by Neural Network Operators}. Vol. 608. Cham: Springer International Publishing, 2016. \url{https://doi.org/10.1007/978-3-319-20505-2}.
		
		\bibitem{Fabra2024} 
		Fabra, Arnau, et al. "Approximation of acoustic black holes with finite element mixed formulations and artificial neural network correction terms." \textit{Finite Elements in Analysis and Design} 241 (2024): 104236. \url{https://doi.org/10.1016/j.finel.2024.104236}.
		
	\end{thebibliography}
\end{document}